\documentclass[12pt]{article}
\usepackage{hyperref}
\usepackage{amssymb,amsmath,amsthm}
\usepackage{enumerate}
\usepackage[letterpaper,hmargin=2.70cm,vmargin=3.2cm]{geometry}

\usepackage{setspace}
\makeatletter
\renewcommand\section{\@startsection {section}{1}{\z@}%
                                   {-3.5ex \@plus -1ex \@minus -.2ex}%
                                   {2.3ex \@plus.2ex}%
                                   {\normalfont\large\bfseries}}
\makeatother

\newtheorem{proposition}{Proposition}
\newtheorem{lemma}{Lemma}
\newtheorem{corollary}{Corollary}
\theoremstyle{definition}
\newtheorem{definition}{Definition}
\newtheorem{remark}{Remark}

\newcommand{\abs}[1]{\left|#1\right|}
\newcommand{\norm}[1]{\left\|#1\right\|}
\newcommand{\I}{{\rm i}}
\newcommand{\inner}[2]{\left\langle#1,#2\right\rangle}

\newcommand{\cH}{{\cal H}}
\newcommand{\cc}[1]{\overline{#1}}

\numberwithin{equation}{section}

\DeclareMathOperator{\re}{Re}
\DeclareMathOperator{\im}{Im}
\DeclareMathOperator{\dom}{Dom}
\DeclareMathOperator{\Ker}{Ker}

\DeclareMathOperator*{\res}{Res}

\DeclareMathOperator{\Sp}{Sp}
\DeclareMathOperator{\Span}{Span}
\begin{document}

\title{Bounded Rank-one Perturbations in Sampling Theory%
\thanks{%
Mathematics Subject Classification(2000):
41A05, 
46E22, 
47A55, 
47B25, 
47N99, 
94A20} 
\thanks{%
Keywords:
rank-one perturbations, sampling expansions}
\thanks{%
Research partially supported by CONACYT under Project P42553\-F.}
\\[6mm]}
\author{\textbf{Luis O. Silva and Julio H. Toloza}
\\[6mm]
Departamento de M\'{e}todos Matem\'{a}ticos y Num\'{e}ricos\\
Instituto de Investigaciones en Matem\'aticas Aplicadas y en Sistemas\\
Universidad Nacional Aut\'onoma de M\'exico\\
C.P. 04510, M\'exico D.F.
\\[4mm]
  \texttt{silva@leibniz.iimas.unam.mx}\\
  \texttt{jtoloza@leibniz.iimas.unam.mx}} \date{} \maketitle
\begin{center}
\begin{minipage}{5in}
  \centerline{{\bf Abstract}} \bigskip
  Sampling theory concerns the problem of reconstruction of functions
  from the knowledge of their values at some discrete set of points.
  In this paper we derive an orthogonal sampling theory and
  associated Lagrange interpolation formulae from a family of bounded
  rank-one perturbations of a self-adjoint operator that has only
  discrete spectrum of multiplicity one.
\end{minipage}
\end{center}
\newpage
\section{Introduction}
\label{sec:intro}

Sampling theory is concerned with the problem of reconstruction of
functions, in a pointwise manner, from the knowledge of their values
at a prescribed discrete set of points.  Resolution of concrete
situations in this theory generally involves the characterization of a
class (usually a linear set) of functions to be interpolated, the
specification of a set of sampling points to be used for all the
functions in the given class, and the derivation of an interpolation
formula.

The cornerstone for many works on sampling theory is the Kramer
sampling theorem \cite{kramer} and its analytic extension
\cite{everitt1}.  Orthogonal sampling formulae often arise as
realizations of this theorem.  The celebrated
Whittaker-Shannon-Kotel'nikov sampling theorem
\cite{kotelnikov,shannon,whittaker} is also a particular case of the
Kramer theorem, although historically the former came first and
motivated the latter.

Orthogonal sampling formulae have been obtained in connection with
differential and difference self-adjoint boundary value problems (see
for instance \cite{garcia00,zayed2} and, of course, the paper due to
Kramer himself \cite{kramer}), and also by resorting to Green's
functions methods \cite{annaby1,zayed1}, among other ODE's techniques.
These results suggest that the spectral theory of operators should
provide a unifying approach to sampling theory.  Following this idea,
a general method for obtaining analytic, orthogonal sampling formulae
has been derived in \cite{nosotros} on the basis of the theory of
representation of simple symmetric operators due to M. G. Krein
\cite{krein1,krein2,krein3,krein4}.  Roughly speaking, the technique
given in \cite{nosotros} consists in the following: By
\cite{krein1,krein2,krein3}, every closed simple symmetric operator
$A$ in a Hilbert space $\cH$ generates a bijective isomorphism between
$\cH$ and a space of functions $\widehat{\cH}$ with certain analytic
properties.  If the operator $A$ satisfies some additional conditions,
all of its self-adjoint extensions have discrete spectrum and every
function $f$ in $\widehat{\cH}$ can be uniquely reconstructed, as long
as one knows the value of $f$ at the spectrum of any self-adjoint
extension of $A$.

Loosely speaking, the self-adjoint extensions of a symmetric operator
with deficiency indices $(1,1)$ constitute a family of singular
rank-one perturbations of one of these self-adjoint extensions
\cite[Sec. 1.1--1.3]{kurasov}.  Therefore, the methods developed in
\cite{nosotros} also holds for singular rank-one perturbations,
provided that these operators correspond to a family of self-adjoint
extensions of some simple symmetric (hence densely defined) operator.
We cannot use, however, a family of bounded rank-one perturbations in
applications to sampling theory without making substantial changes to
the results of \cite{nosotros}.  Krein's approach to symmetric
operators with equal deficiency indices does not work for bounded
rank-one perturbations since operators of this kind may only be seen
as self-adjoint extensions of a certain not densely defined, Hermitian
operator \cite[Sec. 1.1]{kurasov}.

The main motivation of the present work is to develop a method in
sampling theory analogous to \cite{nosotros} for the case of bounded
rank-one perturbations.  With this purpose in mind, we begin by
constructing a representation space for bounded rank-one perturbations
in some sense similar to that of Krein for simple symmetric operators.
Elements of this representation space are the functions to be
interpolated.  Then we obtain a general Kramer-type analytic sampling
formula which turns out to be a Lagrange interpolation formula.  We
also characterize the space of interpolated functions as a space of
meromorphic functions with some properties resembling those of a de
Brange space.  Examples are discussed in the last part of this work.

\section{Preliminaries}
\label{sec:prel}

The following review is based on the standard treatment concerning
rank-one perturbations of self-adjoint operators, as discussed in
detail by Gesztesy and Simon \cite{gesztesy, simon1}, and
Albeverio and Kurasov \cite{kurasov}.

In a separable Hilbert space $\cH$, we consider a possibly unbounded,
self-adjoint operator $A$ with discrete spectrum of multiplicity one.
Let $\mu$ be a cyclic vector for $A$, that is,
$\{(A-zI)^{-1}\mu:z\in\mathbb{C}\}$ is a total set in
$\cH$. Throughout this work we assume that $\norm{\mu}=1$.

Given $\mu$, let us define the family of bounded rank-one perturbations
of $A$,
\begin{equation}
  \label{eq:family}
  A_h:=A+h\langle\mu,\cdot\rangle\mu\,,\qquad h\in\mathbb{R}\,,
\end{equation}
where the inner product is taken, from now on, anti-linear in its
first argument.  Naturally, $\dom (A_h)=\dom (A)$ for any
$h\in\mathbb{R}$.  Elementary perturbation theory implies that all the
elements of (\ref{eq:family}) have discrete spectrum.  As pointed out
in \cite[Sec. 1.1]{kurasov}, the operators $A_h$ may be seen as
self-adjoint extensions of some Hermitian operator, in the sense of
\cite{weidmann}, with non-dense domain.

Consider the family of functions
\begin{equation}
  \label{eq:weyl-function}
  F_h(z):=\langle\mu,(A_h-zI)^{-1}\mu\rangle\,,\quad
  z\not\in\Sp(A_h)\,,\quad h\in\mathbb{R}\,.
\end{equation}
In the sequel we shall denote $F_0$ by $F$. By the spectral theorem,
$F_h(z)$ is the Borel transform of the spectral function
$m_h(t)=\langle\mu,E_h(t)\mu\rangle$, where $E_h(t)$ is the
spectral resolution of the identity corresponding to
$A_h$. Hence, $F_h(z)$ is a Herglotz meromorphic function having simple
poles at the eigenvalues of $A_h$.

From the second resolvent identity \cite[Thm. 5.13]{weidmann} one obtains
\begin{equation*}
   (A_h-zI)^{-1}= (A-zI)^{-1}-h
  \langle(A_h-\cc{z}I)^{-1}\mu,\cdot\rangle(A-zI)^{-1},\quad h\in\mathbb{R}\,.
\end{equation*}
This equation yields the well-known Aronzajn-Krein formula
\cite[Eq. 1.3]{simon1}
\begin{equation}
  \label{eq:a-k-formula}
   F_h(z)=\frac{F(z)}{1+hF(z)},\qquad h\in\mathbb{R}\,.
\end{equation}
The spectral properties of the whole family (\ref{eq:family}) are
contained in (\ref{eq:a-k-formula}). Indeed, one can easily show that
the function $F_h/F_{h'}$, $h\ne h'$, is also Herglotz and its zeros
and poles are given by the poles of $F_{h'}$ and $F_h$ respectively
\cite{weder}. Thus, the spectra of any two different elements of the
family (\ref{eq:family}) interlace, i.\,e., between two neighboring
eigenvalues of one operator there is one and only one eigenvalue of
any other. Also, (\ref{eq:a-k-formula}) implies that
$x_0\in\mathbb{R}$ is a pole of $F_h$ if and only if
\begin{equation}
  \label{eq:pole-f-h}
  \frac{1}{F(x_0)} + h=0\,.
\end{equation}
Therefore for any $x\in\mathbb{R}$ which is not a zero of $F$, there
exists a unique $h\in\mathbb{R}$ such that $x$ is an eigenvalue of
$A_h$. One can extend this result to every $x\in\mathbb{R}$ by
considering an infinite coupling constant $h=\infty$ in
(\ref{eq:family}) (see \cite[Sec 1.5]{simon1},
\cite[Sec. 1.1.2]{kurasov}).  From the properties of $F(z)$, it is
shown that $A_\infty$ also have simple discrete spectrum (see footnote
in \cite[p. 55]{akhiezer1}) and
$\Sp(A_\infty)=\{x\in\mathbb{R}:F(x)=0\}$.  Thus, for any
$x\in\mathbb{R}$, there exists a unique $h\in\mathbb{R}\cup\{\infty\}$
such that $x$ is an eigenvalue of $A_h$.

The main peculiarity of $A_\infty$, that separates it from the family
$A_h$ with finite $h$, is its domain. Indeed, for the kind of
perturbations considered here, the domain of $A_\infty$ is the set
$\{\varphi\in\dom(A):\langle\varphi,\mu\rangle=0\}$
\cite[Sec. 1.1.1]{kurasov}, \cite[Thm.~1.15]{simon1}.

\section{Sampling theory}
\label{sec:sampling}
Based on the theory of rank-one perturbations, we construct in this
section a linear space of meromorphic functions $\widehat{\cH}_\mu$
and derive an interpolation formula valid for all the elements in
$\widehat{\cH}_\mu$.

In our considerations below, the following vector-valued function of
complex argument will play an important r\^{o}le.
\begin{equation*}
  \xi(z):=\frac{(A-\cc{z}I)^{-1}\mu}{F(\cc{z})}\,,\qquad
  z\not\in\Sp(A_\infty)\,.
\end{equation*}
Notice that $\xi(z)$ is well defined for $z\in\Sp(A)$.

\begin{lemma}
\label{lem:eigenvectors}
For any $x\in\mathbb{R}\setminus\Sp(A_\infty)$, there exists
(a unique) $h\in\mathbb{R}$ such that
\begin{equation}\label{muh}
	\xi(x)\in\Ker (A_h-xI)\,.
\end{equation}
Similarly,
\begin{equation}\label{muinfty}
	(A-xI)^{-1}\mu\in\Ker (A_\infty-xI)
\end{equation}
for every $x\in\Sp(A_\infty)$.
\end{lemma}
\begin{proof}
  We first consider $x\not\in\Sp(A)\cup\Sp(A_\infty)$. Let $h\ne 0$ be
  such that $x\in\Sp(A_h)$ (we already know that there is always such
  $h$).  We have
\begin{equation*}
\begin{split}
A_h\xi(x)=\frac{1}{F(x)}A_h(A-xI)^{-1}\mu
         =\left(\frac{1}{F(x)} + h\right)\mu+\frac{x}{F(x)}(A-xI)^{-1}\mu\,.
\end{split}
\end{equation*}
The first assertion of the lemma follows from the last expression and
(\ref{eq:pole-f-h}). When $x\in\Sp(A)$, the statement follows by
a limiting argument based on the fact that $A$ is closed.

We now prove the last assertion of the lemma. Define
$P:=\langle\mu,\cdot\rangle\mu$. By virtue of \cite[Thm. 1.18,
Rem. 2]{simon1}, there is a cyclic vector $\eta$ that obeys
\begin{equation}
    \label{eq:a-infinity-action}
    (A_\infty-zI)^{-1}\eta=\frac{1}{F(z)}(I-P)(A-z)^{-1}\mu\,,
    \qquad z\not\in\Sp(A)\cup\Sp(A_\infty)\,.
  \end{equation}
Clearly, $\eta\in\dom (A_\infty)$. We compute the projection of $\eta$ along
the eigenspace associated to $x\in\Sp(A_\infty)$. Using
(\ref{eq:a-infinity-action}) and some $\epsilon>0$ sufficiently small, we
obtain
 \begin{equation*}
   \begin{split}
   \left[E_\infty(x+0)-E_\infty(x-0)\right]\eta
   &=\frac{1}{2\pi\I}\int_{\abs{x-z}=\epsilon}(A_\infty-z)^{-1}\eta\,dz
	\\[1mm]
   &=\int_{\abs{x-z}=\epsilon}\left[\frac{1}{F(z)}(A-zI)^{-1}-I\right]\mu\,dz
	\\[1mm]
   &=\left(\res_{z=x}\frac{1}{F(z)}\right)(A-xI)^{-1}\mu\,.
   \end{split}
 \end{equation*}
Since the last expression is different from zero, (\ref{muinfty}) is proven.\\
\end{proof}

\begin{definition}
  \label{def:space}
  For any $\varphi\in\cH$, let $\Phi_\mu$ be the mapping given by
  \begin{equation*}
    (\Phi_\mu\varphi)(z):=\langle\xi(z),\varphi\rangle\,,\qquad
    z\in\mathbb{C}\setminus\Sp(A_\infty)\,.
  \end{equation*}
We sometimes shall denote $\Phi_\mu\varphi$ by $\widehat{\varphi}$.
\end{definition}

The mapping $\Phi_\mu$ is a linear injective operator from $\cH$ onto
a certain space of meromorphic functions
$\widehat{\cH}_\mu:=\Phi_\mu\cH$. The injectivity may be verified with
the aid of (\ref{muh}). Some properties that characterize
the set $\widehat{\cH}_\mu$ will be accounted for in the next section.

\begin{proposition}
\label{prop:sampling}
 Given some fixed $h\in\mathbb{R}$,
  let $\{x_j\}_j=\Sp(A_h)$. Define $G_h(z):=1/F_h(z)=h+1/F(z)$. Then, for every
  $f(z)\in\widehat{\cH}_\mu$, we have
\begin{equation}
 \label{eq:interpolation}
f(z) = \sum_{x_j\in\Sp(A_h)}
\frac{G_h(z)}{(z-x_j)G_h'(x_j)} f(x_j)\,,
	\qquad z\in\mathbb{C}\setminus\Sp(A_\infty)\,.
\end{equation}
The series is uniformly convergent on every compact subset of the domain.
\end{proposition}

\begin{proof}
  Because of the assertion (\ref{muh}) of Lemma~\ref{lem:eigenvectors},
  $\{\xi(x_j)\}_j$ is a complete orthogonal set in $\cH$. Hence
  \begin{equation}
    \label{eq:kramer-formula}
    \widehat{\varphi}(z)=\inner{\xi(z)}{\varphi}
	=\sum_{x_j\in\Sp(A_h)}
    \frac{\inner{\xi(z)}{\xi(x_j)}}{\norm{\xi(x_j)}^2}\widehat{\varphi}(x_j)\,,
  \end{equation}
where the series converges uniformly on compacts of
$\mathbb{C}\setminus\Sp(A_\infty)$ by virtue of the Cauchy-Schwarz inequality.

Now, the first resolvent identity implies
\[
\inner{\xi(z)}{\xi(w)}=
	(z-\cc{w})^{-1}\left[\frac{1}{F(\cc{w})}-\frac{1}{F(z)}\right].
\]
In conjunction with (\ref{eq:pole-f-h}) and the convention $1/\infty=0$
when $h=0$, the last equation gives rise to the identity
\[
\inner{\xi(z)}{\xi(x_j)}=
	-(z-x_j)^{-1}\left[h+\frac{1}{F(z)}\right].
\]
Finally, notice that
\begin{equation*}
G_h'(w) = -\frac{1}{F^2(w)}F'(w)
        = -\frac{1}{F^2(w)}\inner{\mu}{(A-w)^2\mu}
        = -\inner{\xi(w)}{\xi(\cc{w})}.
\end{equation*}
Evaluation of the last expression at $w=x_j$ yields the desired result.\\
\end{proof}

\begin{remark}
  Equation (\ref{eq:kramer-formula}) is an orthogonal sampling formula
  of Kramer-type \cite{kramer}. Since the function $G_h(z)$ has simple
  zeroes at the points of $\Sp(A_h)$, expression (\ref{eq:interpolation}) 
  is indeed a Lagrange interpolation formula.
\end{remark}

\section{Spaces of interpolated functions}
\label{sec:meromorphic}
The present section is devoted to the characterization of the space of
functions $\widehat{\cH}_\mu$ introduced by means of the mapping
$\Phi_\mu$ of Definition \ref{def:space}. Notice that
$\widehat{\cH}_\mu$ depends on both the operator $A$ and the
cyclic vector $\mu$.

The following statement gives a quite explicit description of
$\widehat{\cH}_\mu$.
\begin{proposition}
  \label{thm:space-descrption}
\[
\widehat{\cH}_\mu =
\left\{f(z)=c+\sum_{x_n\in\Sp(A_\infty)}\frac{c_n}{z-x_n}:\,
       c,c_n\in\mathbb{C},
       \sum_{x_n\in\Sp(A_\infty)}\abs{c_n}^2F'(x_n)<\infty\right\},
\]
where the series above converge uniformly on compact subsets of
$\mathbb{C}\setminus\Sp(A_\infty)$.
\end{proposition}
\begin{proof} Let $G$ denotes the set defined by the right-hand side
of the statement.

Given $x_n\in\Sp(A_\infty)$, it follows from (\ref{muinfty}) of Lemma
\ref{lem:eigenvectors} that $\omega(x_n):=(A-x_nI)^{-1}\mu$ is the
associated eigenvector (up to normalization). Taking into account the
first resolvent identity, we get
\begin{equation*}
\widehat{\omega(x_n)}(z)
	=\frac{1}{F(z)}\langle(A-\cc{z})^{-1}\mu,(A-x_n)^{-1}\mu\rangle
        =\frac{1}{F(z)}\left[\frac{F(z)}{z-x_n}+ F(x_n)\right]
	=\frac{1}{z-x_n}\,.
\end{equation*}
The expression after the first equality above also implies that
$\norm{\omega(x_n)}^2=F'(x_n)$. Recalling the definition of $\dom(A_\infty)$,
it follows that the set
\begin{equation*}
	B:=\{\mu\}\cup
	\left\{\norm{\omega(x_n)}^{-1}\omega(x_n)\right\}_{x_n\in\Sp(A_\infty)}
\end{equation*}
is an orthonormal basis in $\cH$. Now, take an arbitrary element
$\varphi$ of $\cH$ and expand it on the basis $B$. By applying $\Phi_\mu$
to $\varphi$ and using the Cauchy-Schwarz inequality, we conclude that
$\widehat{H}_\mu\subset G$.

The inclusion $G\subset\widehat{H}_\mu$ follows from noticing that any
$f(z)\in G$ is the image under $\Phi_\mu$ of an element in $\cH$
of the form $c\mu + \sum_{x_n\in\Sp(A_\infty)}c_n\omega(x_n)$.\\
\end{proof}

Notice that, by virtue of Proposition \ref{thm:space-descrption}, the
only entire functions in $\widehat{\cH}_\mu$ are the constant
functions. Also, one easily verifies that any constant function in
$\widehat{\cH}_\mu$ is the image under $\Phi_\mu$ of a vector in
$\Span \{\mu\}$.

The following straightforward result shows that the functions in
$\widehat{\cH}_\mu$ share some properties with those in a
de Branges space.
\begin{lemma}
\label{lem:space-properties}
The space $\widehat{\cal H}_\mu$ has the following properties:
\begin{enumerate}
\item[(i).] Assume that $f(z)\in\widehat{\cal H}_\mu$ has a non-real zero $w$.
	Then $g(z):=\frac{z-\cc{w}}{z-w}f(z)$ also belongs to ${\cal H}_\mu$.
\item[(ii).] The evaluation functional $f(\cdot)\mapsto f(z)$ is continuous for
	every $z\in\mathbb{C}\setminus\Sp(A_\infty)$.
\item[(iii).] For every $f(z)\in\widehat{\cal H}_\mu$,
	$g(z):=\cc{f(\cc{z})}$ belongs to ${\cal H}_\mu$.
\end{enumerate}
\end{lemma}
\begin{proof}
  We have $f(z)=\inner{\xi(z)}{\varphi}$ for some
  $\varphi\in\cH$. Given $w$ such that $f(w)=0$, consider
  $\eta=(A-\cc{w}I)(A-wI)^{-1}\varphi$. A short computation yields
  $g(z)=\inner{\xi(z)}{\eta}$, thus showing (i).  Assertion (ii) is
  rather obvious so the proof is omitted. On the basis of
  Proposition~\ref{thm:space-descrption} one
  verifies (iii).
  \\
\end{proof}

In what follows we show that $\widehat{\cH}_\mu$ can be endowed with
several Hilbert space structures, each one determined by the
spectral functions $m_h(x)$, $h\in\mathbb{R}$.

\begin{lemma}
  Let $h\in\mathbb{R}$ and $\{x_j\}_j=\Sp(A_h)$, arranged in non-decreasing
  order. Then the spectral function $m_h(x)$ is given by
  \begin{equation*}
    m_h(x)=\sum_{x_j\le x}\norm{\xi(x_j)}^{-2}\,.
  \end{equation*}
\end{lemma}
\begin{proof}
Let us recall first the following well-known results \cite[Thm. 1.6]{simon1}
\begin{equation*}
\lim_{\epsilon\to 0}\epsilon\re F_h(x+\I\epsilon)= 0,
\end{equation*}
\begin{equation*}
  \lim_{\epsilon\to 0}\epsilon\im F_h(x+\I\epsilon)
  	=\lim_{\epsilon\to 0}
	 \int_{\mathbb{R}}\frac{\epsilon^2\ dm_h(y)}{(y-x)^2+\epsilon^2}
	=m_h(\{x\})\,,
\end{equation*}
and also the identity
\begin{equation*}
(A_h-zI)^{-1}\mu=\frac{1}{1+hF(z)}(A-zI)^{-1}\mu\,.
\end{equation*}

Consider $x_j\in\Sp(A_h)$. It suffices to verify that
$m_h(\{x_j\})=\norm{\xi(x_j)}^{-2}$. By resorting to the equalities mentioned
above, a straightforward computation shows that
\begin{align*}
\inner{\xi(x_j-i\epsilon)}{\xi(x_j-i\epsilon)}
	&= \frac{1}{\abs{F(x_j+i\epsilon)}^2}
	   \inner{(A-(x_j+i\epsilon)I)^{-1}\mu}{(A-(x_j+i\epsilon)I)^{-1}\mu}\\
	&= \frac{1}{\abs{F_h(x_j+i\epsilon)}^2}
	   \inner{(A_h-(x_j+i\epsilon)I)^{-1}\mu}
	         {(A_h-(x_j+i\epsilon)I)^{-1}\mu}\\
	&= \frac{1}{[\epsilon\re F_h(x_j+i\epsilon)]^2
	   +[\epsilon\im F_h(x_j+i\epsilon)]^2}
           \int_{\mathbb{R}}\frac{\epsilon^2\ dm_h(y)}{(y-x_j)^2+\epsilon^2}\\
        &\to\frac{1}{m_h(\{x_j\})}\,,\qquad \epsilon\to 0.
\end{align*}
The proof is now complete.\\
\end{proof}
\begin{remark}
For $h\ne 0$, this result is in fact statement (ii) of
\cite[Thm. 2.2]{simon1} in disguise.
\end{remark}

\begin{proposition}
\label{prop:hilbert-space-structure}
For arbitrary $h\in\mathbb{R}$, the map $\Phi_\mu$ is a unitary transformation
from $\cH$ onto $L^2(\mathbb{R},dm_h)$.
\end{proposition}

\begin{proof}
$\Phi_\mu$ is a linear isometry from $\cH$ into $L^2(\mathbb{R},dm_h)$. Indeed,
\begin{align*}
\inner{\widehat{\varphi}(\cdot)}{\widehat{\psi}(\cdot)}_h
        :=& \int_\mathbb{R}\inner{\varphi}{\xi(x)}\inner{\xi(x)}{\psi}dm_h(x)\\
         =& \sum_{x_j\in\Sp(A_h)}
	    \frac{\inner{\varphi}{\xi(x_j)}\inner{\xi(x_j)}{\psi}}
	    {\norm{\xi(x_j)}^2}
         =  \inner{\varphi}{\psi}.
\end{align*}

Now consider $f(x)\in L^2(\mathbb{R},dm_h)$. This means that
\[
\norm{f(\cdot)}^2_h
	= \sum_{x_j\in\Sp(A_h)}\frac{\abs{f(x_j)}^2}{\norm{\xi(x_j)}^2}
	< \infty.
\]
Define
\[
\eta = \sum_{x_j\in\Sp(A_h)}\frac{f(x_j)}{\norm{\xi(x_j)}^2}\xi(x_j),
\]
which is clearly an element in $\cH$. It is not difficult to verify that
$\norm{f(\cdot)-\widehat{\eta}(\cdot)}_h=0$.\\
\end{proof}

Define $C:=\Phi_\mu^{-1}\widehat{C}\phi_\mu$, where
$(\widehat{C}f)(z)=\cc{f(\cc{z})}$ for $f\in\widehat{\cH}_\mu$. By
(iii) of Lemma~\ref{lem:space-properties} and
Proposition~\ref{prop:hilbert-space-structure}, it follows that $C$ is a
complex conjugation with respect to which both $A$ and $\mu$ are real.

We conclude this section with a comment about the representation of
the operators $A_h$ as operators on $\widehat{\cH}_\mu$.  A simple
computation shows that, for every $h\in\mathbb{R}$, $A_h$ is
transformed by $\Phi_\mu$ into a quasi-multiplication operator, in the
sense that
\begin{equation}\label{mult-op}
\widehat{A_h\varphi}(z)
	= \frac{1}{F_h(z)}\inner{\mu}{\varphi}
	  + z\widehat{\varphi}(z)
\end{equation}
for every $\varphi\in\dom(A)$.  This is obviously the multiplication
operator in $\Phi_\mu\dom(A_\infty)$.  Moreover, (\ref{mult-op}) reduces
to the multiplication operator in a weak sense; indeed,
\[
\inner{\widehat{\varphi}(\cdot)}{\widehat{A_h\psi}(\cdot)}_h
= \inner{\widehat{\varphi}(\cdot)}{(\cdot)\widehat{\psi}(\cdot)}_h\,,
\]
for every $\varphi,\psi\in\dom(A)$.

\section{Examples}
\label{sec:example}
\textbf{Rank-one perturbations of a Jacobi matrix.}
Consider the following semi-infinite Jacobi matrix
\begin{equation}
  \label{eq:jm-0}
  \begin{pmatrix}
    q_1 & b_1 & 0   &  0  &  \cdots \\[1mm]
    b_1 & q_2 & b_2 &  0  &  \cdots \\[1mm]
    0   & b_2 & q_3 & b_3 &  \\
    0 & 0 & b_3 & q_4 & \ddots\\ \vdots & \vdots &  & \ddots & \ddots
  \end{pmatrix}
\end{equation}
with $q_n\in\mathbb{R}$ and $b_n>0$ for $n\in\mathbb{N}$, and define
in the Hilbert space $l^2(\mathbb{N})$ the operator $J$ in such a way
that its matrix representation with respect to the canonical basis
$\{\delta_n\}_{n=1}^\infty$ in $l^2(\mathbb{N})$ is
(\ref{eq:jm-0}). By this definition, $J$
(cf. \cite[Sec.\,47]{akhiezer2}) is the minimal closed symmetric
operator satisfying
\begin{equation*}
  \langle\delta_n,J\delta_n\rangle = q_n\,,
\quad \langle\delta_{n+1},J\delta_n\rangle =
\langle\delta_n,J\delta_{n+1}\rangle = b_n\,,\quad\forall n\in\mathbb{N}\,.
\end{equation*}
The Jacobi operator $J$ may have deficiency indices $(1,1)$ or $(0,0)$
\cite[Chap.\,4 Sec.\,1.2]{akhiezer1}, \cite[Cor.\,2.9]{simon2}. For
this example we consider $J$ to be self-adjoint, i.\,e., the case of
deficiency indices $(0,0)$. We also assume that $J$ has only discrete
spectrum.  Our family of self-adjoint operators is given by
\begin{equation}
  \label{eq:jm-family}
  J_h:=J+h\langle\delta_1,\cdot\rangle\delta_1\,,\qquad h\in\mathbb{R}\,.
\end{equation}
It is relevant to note that $\delta_1$ is a cyclic vector for $J$ since
the matrix elements $b_n$ are always assumed to be different from zero.

One can study $J$ through the following second order difference system
\begin{equation}
  \label{eq:main-recurrence}
 b_{n-1}f_{n-1}+ q_nf_n+
  b_nf_{n+1}=
  z f_n\qquad n>1\,,\ z\in\mathbb{C}\,.
\end{equation}
with boundary condition
\begin{equation}
  \label{eq:boundary}
  q_1f_1+b_1f_2=z  f_1\,.
\end{equation}
If one sets $f_1=1$, then $f_2$ is completely determined by
(\ref{eq:boundary}). Having $f_1$ and $f_2$, the equation
(\ref{eq:main-recurrence}) gives all the other elements of a sequence
$\{f_n\}_{n=1}^\infty$ that formally satisfies
(\ref{eq:main-recurrence}) and (\ref{eq:boundary}). $f_n$ is a
polynomial of $z$ of degree $n-1$, so we denote $f_n=:P_{n-1}(z)$. The
polynomials $P_n(z)$, $n=0,1,2,\dots$ are referred to as the
polynomials of the first kind associated with the matrix
(\ref{eq:jm-0}). The polynomials of the second kind $Q_n(z)$,
$n=0,1,2,\dots$ associated with (\ref{eq:jm-0}) are defined
as the solutions of
\begin{equation*}
    b_{n-1}f_{n-1} + q_n f_n + b_nf_{n+1} = z f_n
\quad n \in \mathbb{N} \setminus \{1\}\\
\end{equation*}
under the assumption that $f_1=0$ and $f_2=b_1^{-1}$. Then
\begin{equation*}
 Q_{n-1}(z):=f_n\,,\quad\forall n\in\mathbb{N}\,.
\end{equation*}
$Q_n(z)$ is a polynomial of degree $n-1$.

Let $P(z)=\{P_n(z)\}_{n=0}^\infty$ and
$Q(z)=\{Q_n(z)\}_{n=0}^\infty$. Then, classical results in the theory
of Jacobi matrices \cite{akhiezer1} give us the following expression
for $\xi(z)$ defined in Section~\ref{sec:sampling}:
\begin{equation*}
  \xi(z)=P(z)+\frac{1}{F(z)}Q(z)\,,
\end{equation*}
where $F(z)$ is the function given by (\ref{eq:weyl-function}) with
$h=0$. In this context $F(z)$ is referred to as the Weyl function of
$J$ and may be determined by
\begin{equation}
  \label{eq:weyl-limit}
  F(z)=-\lim_{n\to\infty}\frac{1}{w_n(z)}\,,\qquad
  w_n(z):=\frac{P_n(z)}{Q_n(z)}\,,
\end{equation}
where the convergence is uniform on any compact subset of
$\mathbb{C}\setminus\Sp(J)$ \cite[Secs.\,2.4,\,4.2]{akhiezer1}.

The operator $J_\infty$ corresponds in this case to the operator in
$l^2(2,\infty)$ whose matrix representation is (\ref{eq:jm-0}) with
the first column and row removed.

For any $f\in\widehat{\cH}_{\delta_1}$ there is a sequence
$\{\varphi_k\}_{k=1}^\infty\in l^2(\mathbb{N})$ such that
\begin{equation*}
  f(z)=\sum_{k=1}^\infty\left(P_{k-1}(z)\varphi_k+
  \frac{\varphi_k}{F(z)}Q_{k-1}(z)\right)\,.
\end{equation*}
Notice that the poles of $f$ are the eigenvalues of $J_\infty$.

By Proposition~\ref{prop:sampling} we have the following interpolation
formula
\begin{equation}
  \label{eq:jm-interpolation}
  f(z)=\lim_{n\to\infty}\sum_{x_j\in\Sp(J_h)}
\frac{h-w_n(z)}{(x_j-z)w_n'(x_j)}f(x_j)\,,\qquad h\in\mathbb{R}\,.
\end{equation}
Indeed, one can write
(\ref{eq:jm-interpolation}) on the basis of (\ref{eq:interpolation})
using the uniform convergence of the limit and the series in (\ref{eq:weyl-limit})
and (\ref{eq:interpolation}), respectively, and the fact that
$w_n'(z)$ is also uniform convergent.
\newline

\noindent
\textbf{One-dimensional harmonic oscillator.} In this example we look
at the sampling formula provided by only the unperturbed operator.

On $L^2(\mathbb{R},dx)$, consider the differential operator
\[
A:=-\frac{d^2}{dx^2} + x^2\,,
\]
which is essentially self-adjoint on $C^\infty_0(\mathbb{R})$. The eigenvalues
are $2n+1$ for $n\in\mathbb{N}\cup\{0\}$; the corresponding eigenfunctions are
\[
\phi_n(x)=\pi^{-1/4}(2^n n!)^{-1/2}e^{-x^2/2}H_n(x),
\]
where $H_n(x)$ are the Hermite polynomials.

A cyclic vector for the operator $A$ is
\begin{equation}\label{mu-function}
\mu(x)=\sum_{n=0}^\infty \frac{1}{(n!)^{1/2}}\phi_n(x)
      =\frac{1}{\pi^{1/4}}e^{-\frac12(x^2-2\sqrt{2}x+1)}\,,
\end{equation}
where the last equality follows from a quick look to the generating
function of the Hermite polynomials.

Now, for every $z\not\in\Sp(A)$,
\begin{equation*}
F(z) = \sum_{n=0}^\infty\frac{1}{n!(2n+1-z)}\,.
\end{equation*}
An elementary argument involving series of partial fractions then shows that
\begin{equation}\label{weyl-function}
F(z) = \frac{1}{2\cos\frac{\pi z}{2}}
       \int_{-\pi}^\pi e^{-\cos\theta+i\sin\theta}
                       e^{i\frac{1-z}{2}\theta}d\theta\,.
\end{equation}
(See, for instance, \cite{barry}.)

The mapping $\Phi_\mu$ is defined by the vector-valued function
\[
\xi(z;x) = \frac{[(A-\cc{z}I)^{-1}\mu](x)}{F(\cc{z})}
         = \frac{1}{F(\cc{z})}\int_{-\infty}^\infty K(z;x,y)\mu(y) dy\,,
\]
where the Green's function $K(z;x,y)$ is given by (see \cite{titchmarsh})
\begin{equation}\label{green-function}
K(z;x,y) = -\frac{\pi^{1/2}}
           {2\Gamma\left(\frac{1+z}{2}\right)\cos\frac{\pi z}{2}}\times
	   \left\{\begin{array}{ll}
	          D_{\frac{z-1}{2}}(2^{1/2}x)D_{\frac{z-1}{2}}(-2^{1/2}y)\,,&
						y\leq x\,,\\[2mm]
		  D_{\frac{z-1}{2}}(-2^{1/2}x)D_{\frac{z-1}{2}}(2^{1/2}y)\,,&
						y > x\,.
	          \end{array}\right.
\end{equation}
In the last expression, $D_p(x)$ denotes the parabolic cylinder function of
order $p$.

In a fashion more customary for sampling theory, we state the following result:

\begin{corollary}
Let $\mu(x)$, $F(z)$ and $K(z;x,y)$ be given by (\ref{mu-function}),
(\ref{weyl-function}) and (\ref{green-function}), respectively.
Let $G(z):=1/F(z)$. Then, every function $f(z)$ of the form
\[
f(z)=\frac{1}{F(z)}\int_{-\infty}^\infty
     \int_{-\infty}^\infty K(z;x,y)\mu(y)\varphi(x) dydx\,,
\]
with $\varphi(x)\in L^2(\mathbb{R},dx)$, is uniquely determined by
the formula
\[
f(z) = \sum_{n=0}^\infty \frac{G(z)}{(z-2n-1)G'(2n+1)}f(2n+1)\,.
\]
\end{corollary}

\end{document}